# BAYESIAN ANALYSIS FOR REVERSIBLE MARKOV CHAINS

By Persi Diaconis and Silke W. W. Rolles

*Stanford University and Eindhoven University of Technology*

We introduce a natural conjugate prior for the transition matrix of a reversible Markov chain. This allows estimation and testing. The prior arises from random walk with reinforcement in the same way the Dirichlet prior arises from Pólya's urn. We give closed form normalizing constants, a simple method of simulation from the posterior and a characterization along the lines of W. E. Johnson's characterization of the Dirichlet prior.

**1. Introduction.** Modeling with Markov chains is an important part of time series analysis, genomics and many other applications. *Reversible* Markov chains are a mainstay of computational statistics through the Gibbs sampler, Metropolis algorithm and their many variants. Reversible chains are widely used natural models in physics and chemistry where reversibility (often called detailed balance) is a stochastic analog of the time reversibility of Newtonian mechanics.

This paper develops tools for a Bayesian analysis of the transition probabilities, stationary distribution and future prediction of a reversible Markov chain. We observe $X_0 = v_0$, $X_1 = v_1, \dots, X_n = v_n$ from a reversible Markov chain with a finite state space $V$. Neither the stationary distribution $\nu(v)$ nor the transition kernel $k(v, v')$ is assumed known. Reversibility entails $\nu(v)k(v, v') = \nu(v')k(v', v)$ for all $v, v' \in V$. We also assume we know which transitions are possible [for which $v, v' \in V$ is $k(v, v') > 0$].

In Section 2 we introduce a family of natural conjugate priors. These are defined via closed form densities and by a generalization of Pólya's urn to random walk with reinforcement on a graph. The density gives normalizing constants needed for testing independence versus reversibility or reversibility versus a full Markovian specification. The random walk gives a simple method of simulating from the posterior (Section 4.5).

Properties of the prior are developed in Section 4. The family is closed under sampling (Proposition 4.1). Mixtures of our conjugates are shown to be dense (Proposition 4.5). A characterization of the priors via predictive properties of the posterior is given (Section 4.2).

A practical example is given in Section 5. Several simple hypotheses are tested for a data set arising from the DNA of the human HLA-B gene. Section 5 also contains remarks about statistical analysis for reversible chains.

**2. A class of prior distributions.** We observe $X_0 = v_0$, $X_1 = v_1, \ldots,$ $X_n = v_n$ from a reversible Markov chain with a finite state space $V$ and unknown transition kernel $k(\cdot,\cdot)$.

Let $G = (V, E)$ be the finite graph with vertex set $V$ and edge set $E$ defined as follows: $e = \{v, v'\} \in E$ (i.e., there is an edge between $v$ and $v'$) if and only if $k(v, v') > 0$. We assume that $k(v, v') > 0$ iff $k(v', v) > 0$. In particular, all edges of $G$ are undirected and an edge is denoted by the set of its endpoints. For some vertices $v$, we may have $k(v, v) > 0$. Define the simplex

$$\Delta := \left\{ x = (x_e)_{e \in E} \in (0,1]^E : \sum_{e \in E} x_e = 1 \right\}. \tag{1}$$

Remark 2.1. The distribution of a reversible Markov chain can be described by putting on the edge between $v$ and $v'$ the weight $x_{\{v,v'\}} := \nu(v)k(v, v') = \nu(v')k(v', v)$. If the weights are normalized so that $\sum_{e \in E} x_e = 1$, this is a unique way to describe the distribution of the Markov chain. A transition from $v$ to $v'$ is made with probability proportional to the weight $x_{\{v,v'\}}$.

Denote by $Q_{v_0,x}$ the distribution of the Markov chain induced by the weights $x = (x_e)_{e \in E} \in \Delta$ which starts with probability 1 in $v_0$. Using this notation, our assumption says that the observed data comes from a distribution in the class

$$\mathcal{Q} := \{Q_{v_0,x} : v_0 \in V, x \in \Delta\}. \tag{2}$$

2.1. *A minimal sufficient statistic.* If the endpoints of an edge $e$ agree, we call $e$ a *loop*. Let

$$E_{\text{loop}} := \{e \in E : e \text{ is a loop}\}. \tag{3}$$

For an edge $e$, denote the set of its endpoints by $\bar{e}$. For $x = (x_e)_{e \in E} \in (0, \infty)^E$ and a vertex $v$, define $x_v$ to be the sum of all components $x_e$ with $e$ incident to $v$,

$$x_v := \sum_{\{e : v \in \bar{e}\}} x_e. \tag{4}$$

In sums such as this the sum is over edges including loops.

Let $\pi := (\pi_0, \pi_1, \ldots, \pi_n)$ be an admissible path in $G$. Define

$$k_v(\pi) := |\{i \in \{1, 2, \ldots, n\} : (v, \pi_i) = (\pi_{i-1}, \pi_i)\}| \qquad \text{for } v \in V, \tag{5}$$

$$k_e(\pi) := \begin{cases} |\{i \in \{1, 2, \ldots, n\} : \{\pi_{i-1}, \pi_i\} = e\}|, & \text{for } e \in E \setminus E_{\text{loop}}, \\ 2 \cdot |\{i \in \{1, 2, \ldots, n\} : \{\pi_{i-1}, \pi_i\} = e\}|, & \text{for } e \in E_{\text{loop}}. \end{cases} \tag{6}$$

That is, $k_v(\pi)$ equals the number of times the path $\pi$ leaves vertex $v$; for an edge $e$ which is not a loop, $k_e(\pi)$ is the number of traversals of $e$ by $\pi$, and for a loop $e$, $k_e(\pi)$ is twice the number of traversals of $e$. Recall that the edges are undirected; hence, $k_e(\pi)$ counts the traversals of $e$ in *both* directions. Set

$$Z_n := (X_0, X_1, \ldots, X_n). \tag{7}$$

PROPOSITION 2.2. *The vector of transition counts $(k_e(Z_n))_{e \in E}$ is a minimal sufficient statistic for the model $\mathcal{Q}_{v_0} := \{Q_{v_0,x} : x \in \Delta\}$.*

PROOF. Let $\pi$ be an admissible path in $G$. In order to prove that $(k_e(Z_n))_{e \in E}$ is a sufficient statistic, we need to show that

$$Q_{v_0,x}(Z_n = \pi | (k_e(Z_n))_{e \in E}) \tag{8}$$

does not depend on $x$. If $\pi$ does not start in $v_0$, (8) equals zero. Otherwise, we have

$$Q_{v_0,x}(Z_n = \pi) = \frac{\prod_{e \in E \setminus E_{\text{loop}}} x_e^{k_e(\pi)} \prod_{e \in E_{\text{loop}}} x_e^{k_e(\pi)/2}}{\prod_{v \in V} x_v^{k_v(\pi)}}. \tag{9}$$

It is not hard to see that $k_v(\pi)$ can be expressed in terms of the $k_e(\pi)$ and the first observation $v_0$. Hence, the $Q_{v_0,x}$-probability of $\pi$ depends only on $k_e(\pi)$, $e \in E$, and $v_0$. Thus, (8) equals one divided by the number of admissible paths $\pi'$ with starting point $v_0$ and $k_e(\pi') = k_e(\pi)$ for all $e \in E$, which is independent of $x$.

Suppose $K := (k_e)_{e \in E}$ is not minimal. Then there exists a sufficient statistic $K'$ which needs less information than $K$. Consequently, there exist two admissible paths $\pi$ and $\pi'$ starting in $v_0$ such that $K(\pi) \neq K(\pi')$ and $K'(\pi) = K'(\pi')$. Then

$$\begin{aligned} &\frac{Q_{v_0,x}(Z_n = \pi | K'(Z_n) = K'(\pi))}{Q_{v_0,x}(Z_n = \pi' | K'(Z_n) = K'(\pi'))} \\ &\quad = \frac{Q_{v_0,x}(Z_n = \pi)}{Q_{v_0,x}(Z_n = \pi')} \qquad (10) \\ &\quad = \prod_{e \in E \setminus E_{\text{loop}}} x_e^{k_e(\pi) - k_e(\pi')} \prod_{e \in E_{\text{loop}}} x_e^{(k_e(\pi) - k_e(\pi'))/2} \prod_{v \in V} x_v^{k_v(\pi') - k_v(\pi)}. \end{aligned}$$

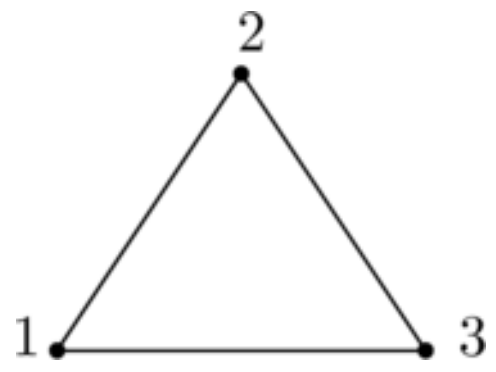

FIG. 1. *The triangle.*

Since by assumption $(k_e(\pi))_{e\in E} \neq (k_e(\pi'))_{e\in E}$, the last quantity depends on $x$. This contradicts the fact that $K'$ is a sufficient statistic. □

2.2. *Definition of the prior densities.* Our aim is to define a class of prior distributions in terms of measures on $\Delta$. We prepare the definition with some notation. We illustrate the definitions by considering a three state process with states $\{1,2,3\}$, all transitions possible, but no holding. This leads to the graph in Figure 1.

Denote the cardinality of a set $S$ by $|S|$. Recall the definition (3) of the set $E_{\text{loop}}$. Set

$$l := |V| + |E_{\text{loop}}| \quad \text{and} \quad m := |E|. \tag{11}$$

For the three state example, $m = l = 3$.

REMARK 2.3. There is a simple way to delineate a *generating set of cycles* of $G$. We call a maximal subgraph of $G$ which contains all loops but no cycle a *spanning tree* of $G$. Choose a spanning tree $T$. Each edge $e \in E \setminus E_{\text{loop}}$ which is not in $T$ forms a cycle $c_e$ when added to $T$. (By definition, a loop is never a cycle and never contained in a cycle.) There are $m-l+1$ such cycles and we enumerate them arbitrarily: $c_1, \dots, c_{m-l+1}$. This set of cycles forms an additive basis for the homology $H_1$ and also serves for our purposes.

For the three state example, we may choose $T$ to have edges $\{1,2\}$ and $\{1,3\}$. Then there is one cycle $c_1$ oriented (say) $1 \to 2 \to 3 \to 1$. In Section 3.4 we show how such a basis of cycles can be obtained for the complete graph.

In general, the first Betti number $\beta_1$ is the dimension of $H_1$. For the complete graph, $\beta_1(K_n) = \binom{n-1}{2}$. Further details can be found in [8], Section 1.16.

DEFINITION 2.4. Orient the cycles $c_1, \dots, c_{m-l+1}$ and all edges $e \in E$ in an arbitrary way. For every $x \in \Delta$, define a matrix $A(x) = (A_{i,j}(x))_{1\le i,j\le m-l+1}$ by

$$A_{i,i}(x) = \sum_{e\in c_i} \frac{1}{x_e}, \qquad A_{i,j}(x) = \sum_{e\in c_i\cap c_j} \pm\frac{1}{x_e} \qquad \text{for } i \neq j, \tag{12}$$

where the signs in the last sum are chosen to be $+1$ or $-1$ depending on whether the edge $e$ has in $c_i$ and $c_j$ the same orientation or not.

In the three state example, the matrix $A(x)$ is $1 \times 1$ with entry $x_{\{1,2\}}^{-1} + x_{\{2,3\}}^{-1} + x_{\{1,3\}}^{-1}$.

Recall the definition (4) of $x_v$. Similarly, define $a_v$ for $a := (a_e)_{e \in E} \in (0,\infty)^E$. The main definition of this section (the conjugate prior) follows.

Definition 2.5. For all $v_0 \in V$ and $a := (a_e)_{e\in E} \in (0,\infty)^E$, define

$$\phi_{v_0,a}(x) := Z_{v_0,a}^{-1} \frac{\prod_{e\in E\setminus E_{\text{loop}}} x_e^{a_e - 1/2} \prod_{e \in E_{\text{loop}}} x_e^{(a_e/2)-1}}{x_{v_0}^{a_{v_0}/2} \prod_{v \in V\setminus\{v_0\}} x_v^{(a_v+1)/2}} \sqrt{\det(A(x))} \tag{13}$$

for $x := (x_e)_{e\in E} \in \Delta$ with

$$\begin{aligned} Z_{v_0,a} := {} & \frac{\prod_{e\in E} \Gamma(a_e)}{\Gamma(a_{v_0}/2) \prod_{v\in V\setminus\{v_0\}} \Gamma((a_v+1)/2) \prod_{e\in E_{\text{loop}}} \Gamma((a_e+1)/2)} \\ & \times \frac{(m-1)!\pi^{(l-1)/2}}{2^{1-l+\sum_{e\in E} a_e}}. \end{aligned} \tag{14}$$

For the three state example with parameters $x_{\{1,2\}} = x$, $x_{\{2,3\}} = y$, $x_{\{1,3\}} = z$, $a_{\{1,2\}} = a$, $a_{\{2,3\}} = b$, $a_{\{1,3\}} = c$, and $v_0 = 1$, the prior density $\phi$ is

$$Z^{-1} \frac{x^{a-1/2} y^{b-1/2} z^{c-1/2}}{(x+z)^{(a+c)/2} (x+y)^{(a+b+1)/2} (y+z)^{(b+c+1)/2}} \sqrt{\frac{1}{x} + \frac{1}{y} + \frac{1}{z}}, \tag{15}$$

with the normalizing constant $Z$ given by

$$\frac{\Gamma(a)\Gamma(b)\Gamma(c)}{\Gamma((a+c)/2)\Gamma((a+b+1)/2)\Gamma((b+c+1)/2)} \cdot \frac{2\pi}{2^{a+b+c-2}}. \tag{16}$$

A derivation of the formula for the density in this special case can be found, for example, in [12]. The density for the triangle with loops is given in (31).

The following proposition shows that the definition of $\phi_{v_0,a}$ is independent of the choice of cycles $c_i$ used in the definition of $A(x)$.

Proposition 2.6. *For the matrix $A$ of Definition 2.4, with $\mathcal{T}$ the set of spanning trees of $G$,*

$$\det A(x) = \sum_{T \in \mathcal{T}} \prod_{e \notin E(T)} \frac{1}{x_e}. \tag{17}$$

PROOF. This identity is proved for graphs without loops in [14], page 145, Theorem 3′. By definition, $A(x)$ does not depend on $x_e$, $e \in E_{\text{loop}}$. Furthermore, since every spanning tree contains all loops, the right-hand side of (17) does not depend on $x_e$, $e \in E_{\text{loop}}$ either. In particular, both sides of (17) are the same for $G$ and the graph obtained from $G$ by removing all loops; hence, they are equal. $\square$

The prior density $\phi_{v_0,a}$ arises in a natural extension of Pólya's urn. We treat this topic next.

2.3. *Random walk with reinforcement.* Let $\sigma$ denote the Lebesgue measure on $\Delta$, normalized such that $\sigma(\Delta) = 1$. The measures $\phi_{v_0,a}\, d\sigma$ on $\Delta$ arise in the study of edge-reinforced random walk, as was observed by Coppersmith and Diaconis; see [3]. Let us explain this connection:

DEFINITION 2.7. All edges of $G$ are given a strictly positive weight; at time 0 edge $e$ has weight $a_e > 0$. An edge-reinforced random walk on $G$ with starting point $v_0$ is defined as follows: The process starts at $v_0$ at time 0. In each step, the random walker traverses an edge with probability proportional to its weight. Each time an edge $e \in E \setminus E_{\text{loop}}$ is traversed, its weight is increased by 1. Each time a loop $e \in E_{\text{loop}}$ is traversed, its weight is increased by 2.

Denote the set of nonnegative integers by $\mathbb{N}_0$. Let $\Omega$ be the set of all $(v_i)_{i\in\mathbb{N}_0} \in V^{\mathbb{N}_0}$ such that $\{v_i, v_{i+1}\} \in E$ for all $i \in \mathbb{N}_0$. Let $X_n : V^{\mathbb{N}_0} \to V$ denote the projection onto the $n$th coordinate. Recall that $Z_n = (X_0, X_1, \ldots, X_n)$. Denote by $P_{v_0,a}$ the distribution on $\Omega$ of an edge-reinforced random walk with starting point $v_0$ and initial edge weights $a = (a_e)_{e\in E}$.

REMARK 2.8. Let $\alpha_e(Z_n) := k_e(Z_n)/n$ be the proportion of traversals of edge $e$ up to time $n$. For a finite graph without loops, it was observed by Coppersmith and Diaconis that $\alpha(Z_n) := (\alpha_e(Z_n))_{e\in E}$ converges almost surely to a random variable with distribution $\phi_{v_0,a}\, d\sigma$; see [3] and also [13]. In particular, $\phi_{v_0,a}\, d\sigma$ is a probability measure on $\Delta$. This fact is not at all obvious from the definition of $\phi_{v_0,a}$.

It turns out that an edge-reinforced random walk on $G$ is a mixture of reversible Markov chains, where the mixing measure described as a measure on edge weights $(x_e)_{e\in E}$ is given by $\phi_{v_0,a}\, d\sigma$. This is made precise by the following theorem.

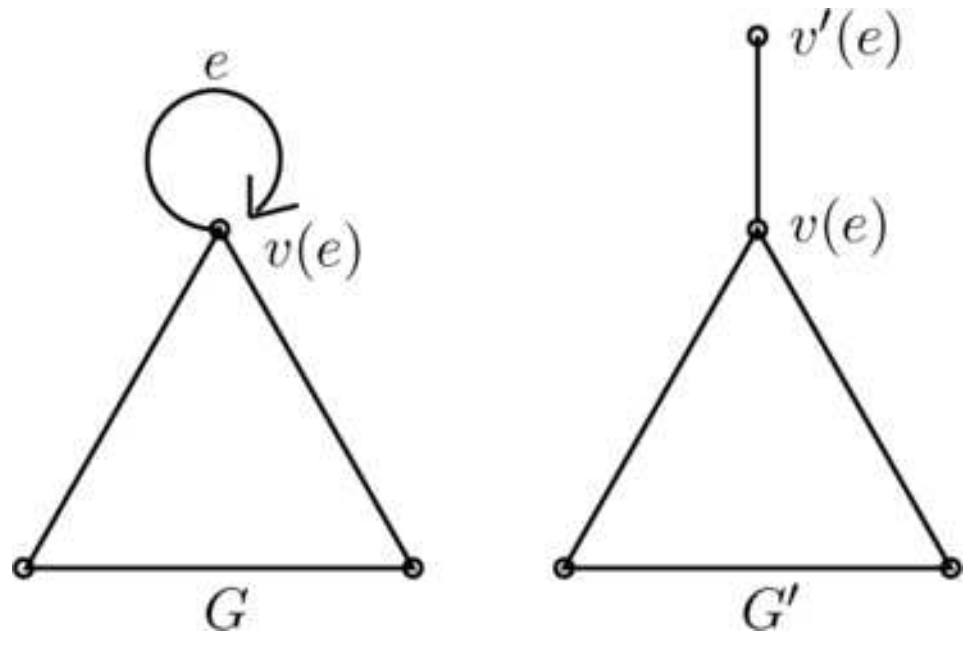

FIG. 2. *Transformation of loops.*

THEOREM 2.9. *Let $(X_n)_{n\in\mathbb{N}_0}$ be an edge-reinforced random walk with initial weights $a=(a_e)_{e\in E}$ starting at $v_0$, and let $Z_n=(X_0,X_1,\ldots,X_n)$. For any admissible path $\pi=(v_0,\ldots,v_n)$, the following holds:*

$$P_{v_0,a}(Z_n=\pi)=\int_\Delta \prod_{i=1}^{n}\frac{x_{\{v_{i-1},v_i\}}}{x_{v_{i-1}}}\phi_{v_0,a}(x)\,d\sigma(x); \tag{18}$$

*here $x:=(x_e)_{e\in E}$. Hence, if $\mathbb{Q}_{v_0,a}$ is the mixture of Markov chains where the mixing measure, described as a measure on edge weights $(x_e)_{e\in E}$, is given by $\phi_{v_0,a}\,d\sigma$, then*

$$P_{v_0,a}=\mathbb{Q}_{v_0,a}. \tag{19}$$

PROOF. If $G$ has no loops, then the claim is true by Theorem 3.1 of [16].

Let $G$ be a graph with loops. Define a graph $G':=(V',E')$ as follows: Replace every loop of $G$ by an edge of degree 1 incident to the same vertex (see Figure 2). More precisely, for all $e\in E_{\mathrm{loop}}$, let $v(e)$ be the vertex $e$ is incident to and let $v'(e)$ be an additional vertex, different from all the others. Then, set $g(e):=\{v(e),v'(e)\}$ and

$$V':=V\cup\{v'(e):e\in E_{\mathrm{loop}}\}, \tag{20}$$

$$E':=[E\setminus E_{\mathrm{loop}}]\cup\{g(e):e\in E_{\mathrm{loop}}\}. \tag{21}$$

The graph $G'$ has no loops and the claim of the theorem is true for $G'$.

Let $P'_{v_0,b}$ be the distribution of a reinforced random walk on $G'$ starting at $v_0$ with initial weights $b=(b_{e'})_{e'\in E'}$ defined by

$$b_{e'}:=\begin{cases} a_{e'}, & \text{if } e'\in E\setminus E_{\mathrm{loop}},\\ a_e, & \text{if } e'=g(e) \text{ for some } e\in E_{\mathrm{loop}}.\end{cases} \tag{22}$$

Any finite admissible path $\pi=(\pi_0=v_0,\pi_1,\ldots,\pi_n)$ in $G$ can be mapped to an admissible path $\pi'=(\pi'_0=v_0,\pi'_1,\ldots,\pi'_{n'})$ in $G'$ by mapping every traversal of a loop $e\in E_{\mathrm{loop}}$ in $\pi$ to a traversal of $(v(e),v'(e),v(e))$ in $\pi'$ [i.e., a traversal of the edge $g(e)$ back and forth in $\pi'$]. The probability that the reinforced random walk on $G$ traverses $\pi$ agrees with the probability that the reinforced random walk on $G'$ traverses $\pi'$. [Note that for $G$ and $G'$ the following is true: Between any two successive visits to $v(e)$, the sum of the weights of all edges incident to $v(e)$ increases by 2.] Since the claim of the theorem is true for $G'$, it follows that

$$\begin{aligned} P_{v_0,a}(Z_n=\pi) &= P'_{v_0,b}(Z_{n'}=\pi') \\ &= \int_\Delta \prod_{i=1}^{n'} \frac{x_{\{\pi'_{i-1},\pi'_i\}}}{x_{\pi'_{i-1}}} \phi'_{v_0,b}(x)\,d\sigma(x), \end{aligned} \tag{23}$$

where $\phi'_{v_0,b}$ denotes the density corresponding to $G'$, starting point $v_0$ and initial weights $b$. We claim that the right-hand side of (23) equals

$$\int_\Delta \prod_{i=1}^{n} \frac{x_{\{\pi_{i-1},\pi_i\}}}{x_{\pi_{i-1}}} \phi_{v_0,a}(x)\,d\sigma(x). \tag{24}$$

Note that a traversal of $e \in E_{\text{loop}}$ contributes $x_e/x_{v(e)}$ to the integrand in (24), whereas a traversal of $(v(e), v'(e), v(e))$ contributes $x_{g(e)}/x_{v(e)}$ to the integrand in (23). Furthermore, $e \in E_{\text{loop}}$ contributes

$$\frac{\Gamma((a_e+1)/2)}{\Gamma(a_e)} \cdot 2^{a_e} \cdot (x_e)^{(a_e/2)-1} \tag{25}$$

to $\phi_{v_0,a}$, whereas the contribution of the edge $g(e)$ and the vertex $v'(e)$ to the density $\phi'_{v_0,b}$ equals

$$\begin{aligned} &\frac{\Gamma((a_{v'(e)}+1)/2)}{\Gamma(a_e)} \cdot 2^{a_e} \cdot \frac{(x_{g(e)})^{a_e-1/2}}{(x_{v'(e)})^{(a_{v'(e)}+1)/2}} \\ &= \frac{\Gamma((a_e+1)/2)}{\Gamma(a_e)} \cdot 2^{a_e} \cdot \frac{(x_{g(e)})^{a_e-1/2}}{(x_{g(e)})^{(a_e+1)/2}} \\ &= \frac{\Gamma((a_e+1)/2)}{\Gamma(a_e)} \cdot 2^{a_e} \cdot (x_{g(e)})^{(a_e/2)-1}. \end{aligned} \tag{26}$$

Finally, $|V| + |E_{\text{loop}}| = |V'|$ and $|E| = |E'|$. Consequently, the expression in (24) agrees with the right-hand side of (23) and the claim follows. □

**3. The prior density for special graphs.** In this section we write down the densities $\phi_{v_0,a}$ for some special graphs.

3.1. *The line graph* (*Birth and death chains*). Consider the line graph with vertex set $V = \{i : 0 \le i \le n\}$ and edge set $E = \{\{i, i+1\} : 0 \le i \le n-1\}$; see Figure 3. Given $a = (a_{\{i-1,i\}})_{1\le i\le n}$, let $b_i := a_{\{i-1,i\}}$. The variables in the simplex $\Delta$ are denoted $z_i := x_{\{i-1,i\}}$.

Recall that the density of the beta distribution with parameters $b_1, b_2 > 0$ is given by

$$\beta[b_1, b_2](p) := \frac{\Gamma(b_1+b_2)}{\Gamma(b_1)\Gamma(b_2)} p^{b_1-1}(1-p)^{b_2-1} \qquad (0 < p < 1). \tag{27}$$

Set

$$p_i = \frac{z_i}{z_i + z_{i+1}}, \qquad 1 \le i \le n-1, \tag{28}$$

and $p := (p_i)_{1\le i\le n-1}$. Clearly, $p_i$ is the probability that the Markov chain with edge weights $z_i$ makes a transition to $i-1$ given it is at $i$. If we make the change of variables (28) in the density $\phi_{v_0,a}$, then we obtain the transformed density $\widetilde{\phi}_{v_0,a}(p)$ given by

$$\widetilde{\phi}_{0,a}(p) = \begin{cases} \displaystyle\prod_{i=1}^{n-1} \beta\left[\frac{b_i+1}{2}, \frac{b_{i+1}}{2}\right](p_i), & \text{if } v_0 = 0, \\ \displaystyle\left[\prod_{i=1}^{v_0-1} \beta\left[\frac{b_i}{2}, \frac{b_{i+1}+1}{2}\right](p_i)\right] \beta\left[\frac{b_{v_0}}{2}, \frac{b_{v_0}}{2}\right](p_{v_0}) & \\ \qquad \displaystyle\times \left[\prod_{i=v_0+1}^{n-1} \beta\left[\frac{b_i+1}{2}, \frac{b_{i+1}}{2}\right](p_i)\right], & \text{if } v_0 \in \{1, 2, \ldots, n-1\}, \\ \displaystyle\prod_{i=1}^{n-1} \beta\left[\frac{b_i}{2}, \frac{b_{i+1}+1}{2}\right](p_i), & \text{if } v_0 = n; \end{cases}$$

here the empty product is defined to be 1.

With the change of variables (28), the conjugate prior can be described as a product of independent beta variables with carefully linked parameters. If loops are allowed, the edge weights are independent Dirichlet by a similar argument (see Section 3.2). The next example contains a generalization.

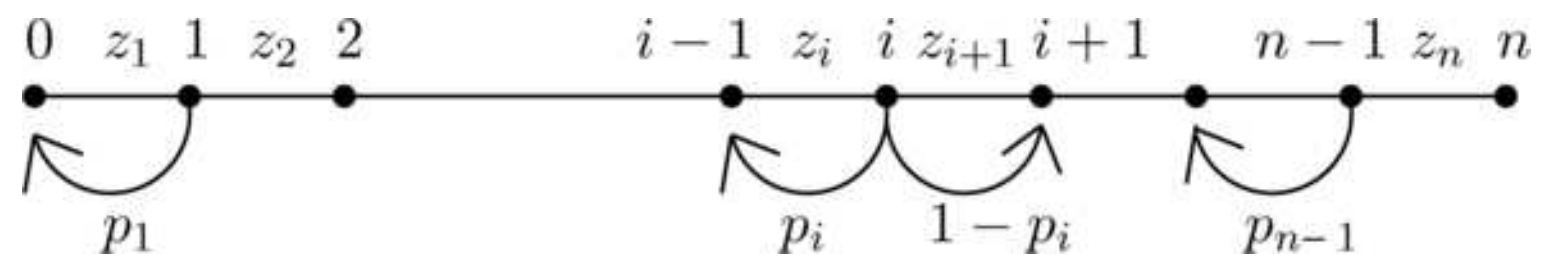

FIG. 3. *The line graph.*

3.2. *Trees with loops.* Recall that the density of the Dirichlet distribution with parameters $b_i > 0$, $1 \le i \le d$, is given by

$$D[b_i; 1 \le i \le d](p_i; 1 \le i \le d) \\ := \frac{\Gamma(\sum_{i=1}^d b_i)}{\prod_{i=1}^d \Gamma(b_i)} \prod_{i=1}^d p_i^{b_i - 1} \qquad \left(p_i > 0, \sum_{i=1}^d p_i = 1\right). \tag{29}$$

Let $T = (V, E)$ be a tree. Suppose that there is a loop attached to every vertex, that is, $\{v\} \in E$ for all $v \in V$. Let $v_0 \in V$. For every $v \in V \setminus \{v_0\}$, there exists a unique shortest path from $v_0$ to $v$. Let $e(v)$ be the unique edge incident to $v$ which is traversed by the shortest path from $v_0$ to $v$. Let $E_v := \{e \in E : v \in \bar{e}\}$ be the set of all edges incident to $v$. Set

$$p_e := \frac{x_e}{x_v} \qquad \text{for } v \in V, e \in E_v, \tag{30}$$

$p := (p_e)_{e \in E}$, and $\vec{p}_v := (p_e)_{e \in E_v}$. If we make the change of variables (30) in the density $\phi_{v_0,a}$, the transformed density $\widetilde{\phi}_{v_0,a}(p)$ is given by

$$D\left[\frac{a_e}{2}, e \in E_{v_0}\right](\vec{p}_{v_0}) \prod_{v \in V \setminus \{v_0\}} D\left[\frac{a_{e(v)} + 1}{2}, \frac{a_e}{2}, e \in E_v \setminus \{e(v)\}\right](\vec{p}_v).$$

Thus again, in the reparametrization (30), the conjugate prior is seen as a product of independent random variables. This is not true in the following example.

The fact that the density $\phi_{v_0,a}$ for a tree has this particular form was first observed by Pemantle [15].

3.3. *The triangle.* Consider the triangle with loops attached to all vertices. Let the vertex set be $V = \{1, 2, 3\}$ and the edge set $E = \{\{1\}, \{2\}, \{3\}, \{1, 2\}, \{1, 3\}, \{2, 3\}\}$ (see Figure 4). Let $b_i$ be the initial weight of the loop at vertex $i$ and let $c_i$ be the initial weight of the edge opposite of vertex $i$. Similarly, let $y_i := x_{\{i\}}$ and let $z_1 := x_{\{2,3\}}$, $z_2 := x_{\{1,3\}}$, $z_3 := x_{\{1,2\}}$.

The density $\phi_{1,a}(y_1, y_2, y_3, z_1, z_2, z_3)$ for $a = (b_1, b_2, b_3, c_1, c_2, c_3)$ is given by

$$Z_{1,a}^{-1} \cdot y_1^{(b_1/2)-1} y_2^{(b_2/2)-1} y_3^{(b_3/2)-1} z_1^{c_1-1} z_2^{c_2-1} z_3^{c_3-1} \sqrt{z_1 z_2 + z_1 z_3 + z_2 z_3} \\ \times \left((y_1 + z_2 + z_3)^{(b_1+c_2+c_3)/2} (y_2 + z_1 + z_3)^{(b_2+c_1+c_3+1)/2} \right. \\ \left. \times (y_3 + z_1 + z_2)^{(b_3+c_1+c_2+1)/2}\right)^{-1}, \tag{31}$$

with

$$Z_{1,a} = \Gamma(c_1)\Gamma(c_2)\Gamma(c_3)\Gamma(b_1/2)\Gamma(b_2/2)\Gamma(b_3/2) \\ \times \left(\Gamma\left(\frac{b_1 + c_2 + c_3}{2}\right)\Gamma\left(\frac{b_2 + c_1 + c_3 + 1}{2}\right) \right. \\ \left. \times \Gamma\left(\frac{b_3 + c_1 + c_2 + 1}{2}\right)\right)^{-1} \cdot \frac{480\pi}{2^{c_1+c_2+c_3}}. \tag{32}$$

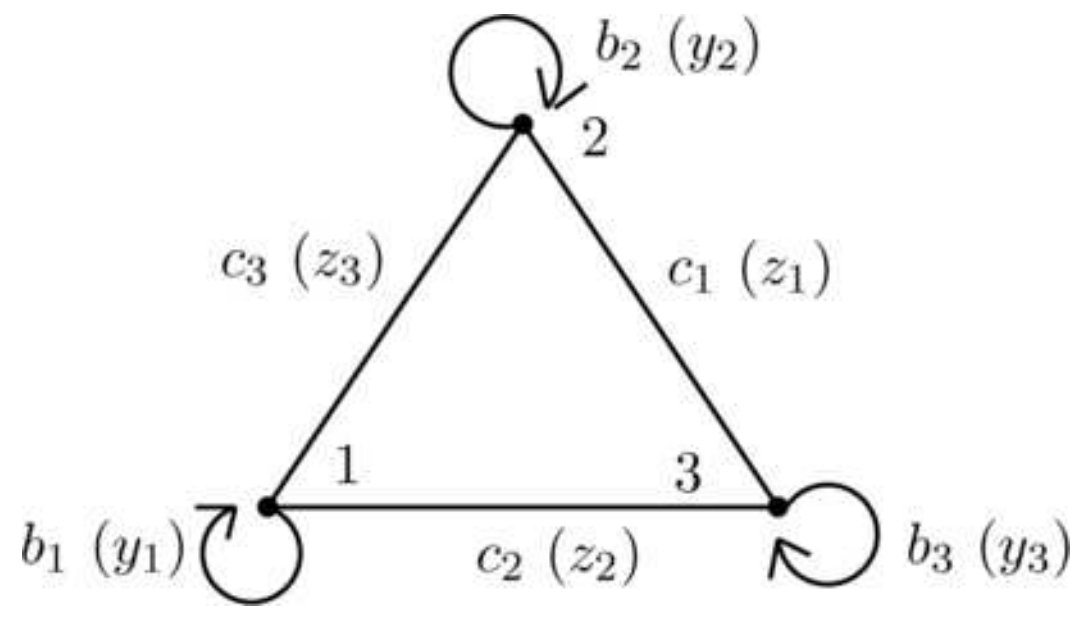

FIG. 4. *The triangle with loops.*

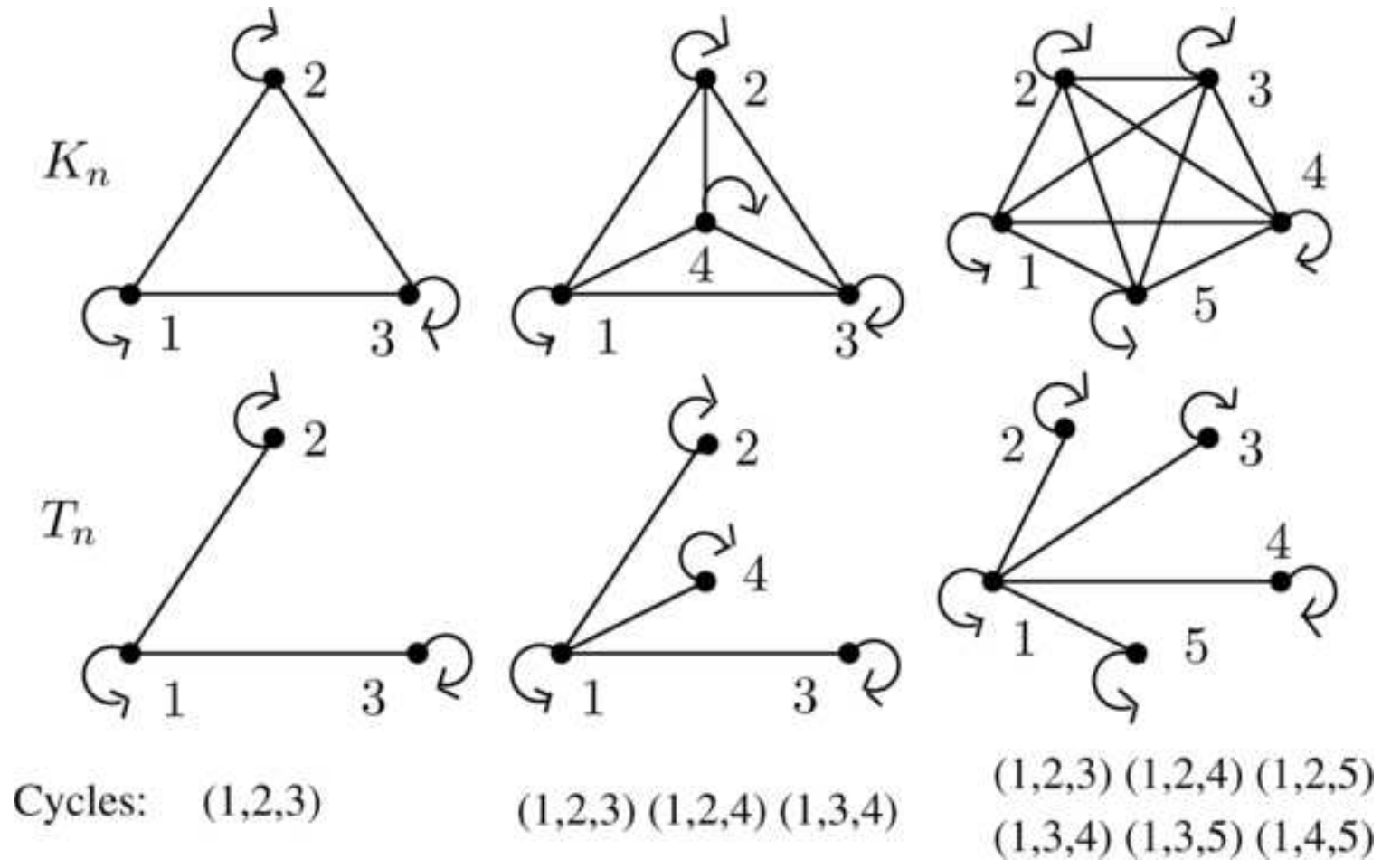

FIG. 5. *The complete graphs $K_3$, $K_4$, $K_5$ with loops and a spanning tree.*

To calculate $Z_{1,a}$ from (14), use the identity

$$\frac{\Gamma(b_i)}{2^{b_i}\Gamma((b_i+1)/2)} = \frac{\Gamma(b_i/2)}{2\sqrt{\pi}} \qquad (i = 1, 2, 3). \tag{33}$$

3.4. *The complete graph.* Perhaps the most important example is where all transitions are possible. This involves the complete graph $K_n$ on $n$ vertices with loops attached to all vertices. Let $V = \{1, 2, 3, \ldots, n\}$. Let $T_n$ be the spanning tree with edges $\{1, i\}$ and loops $\{i\}$, $1 \le i \le n$. This spanning tree induces the basis of cycles given by all triangles $(1, i, j)$, $2 \le i < j \le n$. Figure 5 shows $K_3$, $K_4$ and $K_5$ together with $T_3$, $T_4$ and $T_5$.

We remark that a different basis of cycles is given by $(i, i+1, j)$ for $1 \le i < j + 1 \le n$. This may be proved by induction using the Mayer–Vietoris decomposition theorem based on $K_{n-1}$ and a point.

Let $a = (a_{\{i,j\}})_{1\le i,j\le n}$ be given. For $K_n$, set $b_i := a_{\{i\}}$, $a_i = \sum_{j=1}^n a_{\{i,j\}}$ and $b := \sum_{1\le i,j\le n} a_{\{i,j\}}$. The variables of the simplex are $x = (x_{\{i,j\}})_{1\le i,j\le n}$. Abbreviating $y_i := x_{\{i\}}$ and $x_i = \sum_{j=1}^n x_{\{i,j\}}$, the density $\phi_{1,a}$ is given by

$$\phi_{1,a}(x) = Z_{1,a}^{-1} \cdot \frac{\prod_{1\le i<j\le n} x_{\{i,j\}}^{a_{\{i,j\}}-1/2} \prod_{i=1}^n y_i^{(b_i/2)-1}}{x_1^{a_1/2} \prod_{i=2}^n x_i^{(a_i+1)/2}} \sqrt{\det(A_n(x))}, \tag{34}$$

with $A_n(x)$ defined in (12) and

$$\begin{aligned} Z_{1,a} = {} & \frac{\prod_{1\le i,j\le n} \Gamma(a_{\{i,j\}})}{\Gamma(a_1/2)\prod_{i=2}^n \Gamma((a_i+1)/2)\prod_{i=1}^n \Gamma((b_i+1)/2)} \\ & \times \frac{((n(n+1)/2)-1)!\pi^{n-1/2}}{2^{1-2n+b}}. \end{aligned} \tag{35}$$

**4. Properties of the family of priors.** For $v_0 \in V$ and $a = (a_e)_{e\in E} \in (0,\infty)^E$, abbreviate

$$\mathbb{P}_{v_0,a} := \phi_{v_0,a}\, d\sigma; \tag{36}$$

that is, $\mathbb{P}_{v_0,a}$ is the measure on $\Delta$ with density $\phi_{v_0,a}$. Recall that $\mathbb{Q}_{v_0,a}$ denotes the mixture of Markov chains where the mixing measure, described as a measure on edge weights $(x_e)_{e\in E}$, is given by $\mathbb{P}_{v_0,a}$. In this section we study properties of the set of prior distributions

$$\mathcal{D} := \{\mathbb{P}_{v_0,a} : v_0 \in V, a = (a_e)_{e\in E} \in (0,\infty)^E\}. \tag{37}$$

4.1. *Closure under sampling.* Recall the definition (6) of $k_e(\pi)$ and recall that $Z_n = (X_0, \ldots, X_n)$.

PROPOSITION 4.1. *Under the prior distribution $\mathbb{P}_{v_0,a}$ with observations $X_0 = v_0$, $X_1 = v_1, \ldots,$ $X_n = v_n$, the posterior is given by $\mathbb{P}_{v_n,(a_e+k_e(Z_n))_{e\in E}}$. In particular, the family $\mathcal{D}$ is closed under sampling.*

PROOF. Suppose the prior distribution is $\mathbb{P}_{v_0,a}$ and we are given $n+1$ observations $\pi = (\pi_0, \pi_1, \ldots, \pi_n)$ sampled from $\mathbb{Q}_{v_0,a}$. Then $\pi_0 = v_0$. We claim that the posterior is given by $\mathbb{P}_{\pi_n,(a_e+k_e(\pi))_{e\in E}}$. By Theorem 2.9, $\mathbb{Q}_{v_0,a} = P_{v_0,a}$. The $P_{v_0,a}$-distribution of $\{X_{n+k}\}_{k\ge 0}$ given $Z_n = \pi$ is the distribution of an edge-reinforced random walk starting at the vertex $\pi_n$ with initial values $a_e + k_e(\pi)$. Using the identity (19) again, it follows that the $P_{v_0,a}$-distribution of $\{X_{n+k}\}_{k\ge 0}$ given $Z_n = \pi$ equals $\mathbb{Q}_{\pi_n,(a_e+k_e(\pi))_{e\in E}}$. Thus, the posterior equals $\mathbb{P}_{\pi_n,(a_e+k_e(\pi))_{e\in E}}$, which is an element of $\mathcal{D}$. $\square$

4.2. *Uniqueness.* In this section we give a characterization of our priors along the lines of W. E. Johnson's characterization of the Dirichlet prior. See [18] for history and [19] for a version for nonreversible chains. The closely related topic of de Finetti's theorem for Markov chains is developed by Freedman [7] and Diaconis and Freedman [4]. See also [6].

DEFINITION 4.2. Two finite admissible paths $\pi$ and $\pi'$ are called *equivalent* if they have the same starting point and satisfy $k_e(\pi) = k_e(\pi')$ for all $e \in E$. We define $P$ to be *partially exchangeable* if $P(Z_n = \pi) = P(Z_n = \pi')$ for any equivalent paths $\pi$ and $\pi'$ of length $n$.

For $n \in \mathbb{N}_0$ and $v \in V$, define

$$k_n(v) := |\{i \in \{0, 1, \ldots, n\} : X_i = v\}|. \tag{38}$$

It seems natural to take a class $\mathcal{P}$ of distributions for $(X_n)_{n\in\mathbb{N}_0}$ with the following properties:

P1. For all $P \in \mathcal{P}$, there exists $v_0 \in V$ such that $P(X_0 = v_0) = 1$.
P2. For all $P \in \mathcal{P}$, $v_0$ as in P1, and any admissible path $\pi$ of length $n \geq 1$ starting at $v_0$, we have $P(Z_n = \pi) > 0$.
P3. Every $P \in \mathcal{P}$ is partially exchangeable.
P4. For all $P \in \mathcal{P}$, $v \in V$ and $e \in E$, there exists a function $f_{P,v,e}$ taking values in $[0, 1]$ such that, for all $n \geq 0$,

$$P(X_{n+1} = v|Z_n) = f_{P,X_n,\{X_n,v\}}(k_n(X_n), k_{\{X_n,v\}}(Z_n)).$$

The condition P4 says that, given $X_0, X_1, \ldots, X_n$, the probability that $X_{n+1} = v$ depends only on the following quantities: the observation $X_n$, the number of times $X_n$ has been observed so far, the edge $\{X_n, v\}$ and the number of times transitions between $X_n$ and $v$ (and between $v$ and $X_n$) have been observed so far.

We make the following assumptions on the graph $G$:

G1. For all $v \in V$, $degree(v) \neq 2$.
G2. The graph $G$ is 2-edge-connected, that is, removing an edge does not make $G$ disconnected.

For example, a triangle with loops or the complete graph $K_n$, $n \geq 4$, with or without loops, satisfies G1 and G2, while a path fails both G1 and G2.

Recall that $Q_{v_0,x}$ is the distribution of the reversible Markov chain starting in $v_0$, making a transition from $v$ to $v'$ with probability proportional to $x_{\{v,v'\}}$ whenever $\{v, v'\} \in E$.

THEOREM 4.3. *Suppose the graph $G$ satisfies* G1 *and* G2.

(a) *The set* $\mathcal{M} := \{\mathbb{Q}_{v_0,a} : v_0 \in V, a = (a_e)_{e\in E} \in (0,\infty)^E\}$ *satisfies* P1–P4.

(b) *On the other hand, if* P1–P4 *are satisfied for a set $\mathcal{P}$ of probability distributions, then for all $P \in \mathcal{P}$ there exist $v_0 \in V$ and $a \in (0,\infty)^E$ such that either*

$$\begin{aligned} &P(X_{n+1} = v | Z_n, k_n(X_n) \geq 3) \\ &\qquad = \mathbb{Q}_{v_0,a}(X_{n+1} = v | Z_n, k_n(X_n) \geq 3) \qquad \forall n \geq 0 \quad \textit{or} \\ &P(X_{n+1} = v | Z_n, k_n(X_n) \geq 3) \\ &\qquad = Q_{v_0,a}(X_{n+1} = v | Z_n, k_n(X_n) \geq 3) \qquad \forall n \geq 0. \end{aligned} \tag{39}$$

The second part of the theorem states that either $P$ and $\mathbb{Q}_{v_0,a}$ or $P$ and $Q_{v_0,a}$ essentially agree; only the conditional probabilities to leave from a state which has been visited at most twice could be different.

Proof of Theorem 4.3. It is straightforward to check that $\mathcal{M}$ has the properties P1–P4. For the converse, let $P \in \mathcal{P}$. If $G$ has no loops, then Theorem 1.2 of Rolles [16] implies that there exist $v_0 \in V$ and $a \in (0,\infty)^E$ such that either (39) holds or $P(X_{n+1} = v | Z_n, k_n(X_n) \geq 3) = P_{v_0,a}(X_{n+1} = v | Z_n, k_n(X_n) \geq 3)$ for all $n$. In this case, the claim follows from (19).

If $G$ has loops, consider the graph $G'$ defined in the proof of Theorem 2.9 and the induced process $X' := (X'_n)_{n \in \mathbb{N}_0}$ on $G'$ with reflection at the vertices $v'(e)$, $e \in E_{\text{loop}}$. The process $X'$ satisfies P1–P4. Hence, the claim holds for $X'$ and, consequently, for $(X_n)_{n \in \mathbb{N}_0}$. □

Remark 4.4. The preceding theorem holds under the assumption that the graph $G$ is 2-edge-connected (G2). If $G$ is not 2-edge-connected, a similar statement can be proved for a different class of priors: One replaces the class $\mathcal{D}$ by the mixing measures of a so-called *modified edge-reinforced random walk*; for the definition of this process, see Definition 2.1 of [16]. A uniqueness statement similar to Theorem 4.3 follows from Theorem 2.1 of [16].

4.3. *The priors are dense.* As shown by Dalal and Hall [2] and Diaconis and Ylvisaker [5] for classical exponential families, mixtures of conjugate priors are dense in the space of all priors. This holds for reversible Markov chains.

Proposition 4.5. *The set of convex combinations of priors in $\mathcal{D}$ is weak-star dense in the set of all prior distributions on reversible Markov chains on $G$.*

Proof. For an infinite admissible path $\pi = (\pi_0, \pi_1, \pi_2, \dots)$ in $G$, define $\alpha(\pi) := (\alpha_e(\pi))_{e \in E}$ by $\alpha_e(\pi) := \lim_{n\to\infty} k_e(\pi_0, \pi_1, \dots, \pi_n)/n$ to be the limiting fraction of crossings of the edge $e$ by the path $\pi$. Let $Z_\infty := (X_0, X_1, X_2, \dots)$.

Note that $\alpha(Z_\infty)$ is defined $\mathbb{Q}_{v_0,a}$-a.s. Define $\tau_n$ to be the $n$th return time to $v_0$. Since $G$ is finite, $\tau_n < \infty$ $\mathbb{Q}_{v_0,a}$-a.s. for all $n \in \mathbb{N}$ and all $a \in (0,\infty)^E$.

Let $f : \Delta \to \mathbb{R}$ be bounded and continuous. Denote the expectation with respect to $\mathbb{Q}_{v_0,a}$ by $\mathbb{E}_{v_0,a}$. Since $X_{\tau_n} = v_0$, Theorem 2.9 implies that

$$\mathbb{E}_{v_0,(a_e+k_e(Z_{\tau_n}))_{e\in E}}[f(\alpha(Z_\infty))] = \mathbb{E}_{v_0,a}[f(\alpha(Z_\infty))|Z_{\tau_n}] := M_n. \tag{40}$$

Clearly, $(M_n)_{n\geq 0}$ is a bounded martingale. Hence, by the martingale convergence theorem,

$$\begin{aligned}\lim_{n\to\infty} \mathbb{E}_{v_0,(a_e+k_e(Z_{\tau_n}))_{e\in E}}[f(\alpha(Z_\infty))] &= \mathbb{E}_{v_0,a}[f(\alpha(Z_\infty))|Z_\infty] = f(\alpha(Z_\infty)) \\ &= \int f\, d\delta_{\alpha(Z_\infty)}\end{aligned} \tag{41}$$

$\mathbb{Q}_{v_0,a}$-a.s.; here $\delta_b$ denotes the point mass in $b$. Since $\Delta$ is compact, there is a countable dense subset of the set of bounded continuous functions on $\Delta$. Hence, the above shows that, for $\mathbb{Q}_{v_0,a}$-almost all $Z_\infty$,

$$\mathbb{Q}_{v_0,(a_e+k_e(Z_{\tau_n}))_{e\in E}}(\alpha(Z_\infty) \in \cdot) \Rightarrow \delta_{\alpha(Z_\infty)}(\cdot) \qquad \text{weakly as } n \to \infty. \tag{42}$$

The $\mathbb{Q}_{v_0,a}$-distribution of $\alpha(Z_\infty)$ equals $\mathbb{P}_{v_0,a}$. Thus,

$$\mathbb{P}_{v_0,(a_e+k_e(Z_{\tau_n}))_{e\in E}}(\cdot) \Rightarrow \delta_{\alpha(Z_\infty)}(\cdot) \qquad \text{weakly as } n \to \infty \tag{43}$$

for $\mathbb{Q}_{v_0,a}$-almost all $Z_\infty$. Recall that the $\mathbb{Q}_{v_0,a}$-distribution of $\alpha(Z_\infty)$ (viz., $\mathbb{P}_{v_0,a}$) is absolutely continuous with respect to Lebesgue measure on $\Delta$ with the density $\phi_{v_0,a}$ which is strictly positive in the interior of $\Delta$. Hence, for Lebesgue-almost all $a \in \Delta$, there is a sequence $a_n \in \Delta$ such that $\mathbb{P}_{v_0,a_n} \Rightarrow \delta_a$ weakly. By the Krein–Milman theorem, convex combinations of point masses are weak-star dense in the set of all measures on $\Delta$. Using a standard argument, it follows that the set of convex combinations of the distributions $\mathbb{P}_{v_0,a}$ is dense in the set of all probability measures on $\Delta$. This completes the proof of the proposition. □

4.4. *Computing some moments.* For any edge $e_0 \in E$, we can calculate the probability that the mixture of Markov chains with mixing measure $\phi_{v_0,a}\, d\sigma$ traverses $e_0$ back and forth starting at an endpoint of $e_0$. This gives a closed form for certain moments of the prior $\mathbb{P}_{v_0,a}$.

PROPOSITION 4.6. *For $e_0 \in E \setminus E_{\mathrm{loop}}$ with endpoints $v$ and $v'$, we have*

$$\int_\Delta \frac{(x_{e_0})^2}{x_v x_{v'}} \phi_{v_0,a}(x)\, d\sigma(x) = \begin{cases} \dfrac{a_{e_0}(a_{e_0}+1)}{(a_v+1)(a_{v'}+1)}, & \text{if } v_0 \notin \{v, v'\}, \\ \dfrac{a_{e_0}(a_{e_0}+1)}{a_v(a_{v'}+1)}, & \text{if } v = v_0. \end{cases} \tag{44}$$

*For a loop $e_0 \in E_{\mathrm{loop}}$ incident to $v$, we have*

$$\int_\Delta \frac{x_{e_0}}{x_v} \phi_{v_0,a}(x)\, d\sigma(x) = \begin{cases} \dfrac{a_{e_0}}{a_v+1}, & \textit{if } v \neq v_0, \\ \dfrac{a_{e_0}}{a_v}, & \textit{if } v = v_0. \end{cases} \tag{45}$$

PROOF. *Case* $e_0 \in E \setminus E_{\mathrm{loop}}$. Suppose $v_0$ is an endpoint of $e_0$, say, $v = v_0$. Then

$$\int_\Delta \frac{(x_{e_0})^2}{x_v x_{v'}} \phi_{v_0,a}(x)\, d\sigma(x) = \mathbb{Q}_{v_0,a}(X_0 = v_0, X_1 = v', X_2 = v_0); \tag{46}$$

this is the probability that the mixture of Markov chains traverses the edge $e_0$ back and forth starting at $v_0$. By (19) $\mathbb{Q}_{v_0,a} = P_{v_0,a}$. Hence, (46) equals the probability that an edge-reinforced random walk traverses $e_0$ back and forth, namely,

$$P_{v_0,a}(X_0 = v_0, X_1 = v', X_2 = v_0) = \frac{a_{e_0}(a_{e_0}+1)}{a_v(a_{v'}+1)}. \tag{47}$$

Here we used the fact that the sum of the weights of all edges incident to $v'$ equals $a_{v'} + 1$ after $e_0$ has been traversed once. This proves the claim in the case $v_0 \in \bar{e}_0$.

Suppose $v_0 \notin \bar{e}_0$. Define $b := (b_e)_{e \in E}$ by $b_{e_0} := a_{e_0} + 2$ and $b_e := a_e$ for $e \in E \setminus \{e_0\}$. Then, using the definition of $\phi_{v_0,a}$, we obtain

$$\frac{(x_{e_0})^2}{x_v x_{v'}} \phi_{v_0,a}(x) = \frac{Z_{v_0,b}}{Z_{v_0,a}} \phi_{v_0,b}(x) \qquad \text{for all } x \in \Delta. \tag{48}$$

Using the definition of the normalizing constants $Z_{v_0,a}$ and $Z_{v_0,b}$ and the identity $\Gamma(z+1) = z\Gamma(z)$, it follows that

$$\frac{Z_{v_0,b}}{Z_{v_0,a}} = \frac{\Gamma((a_v+1)/2)\Gamma((a_{v'}+1)/2)\Gamma(a_{e_0}+2)}{4\Gamma((a_v+3)/2)\Gamma((a_{v'}+3)/2)\Gamma(a_{e_0})} = \frac{a_{e_0}(a_{e_0}+1)}{(a_v+1)(a_{v'}+1)}. \tag{49}$$

Since $\int_\Delta \phi_{v_0,b}(x)\, d\sigma(x) = 1$, the claim follows by integrating both sides of (48) over $\Delta$.

*Case* $e_0 \in E_{\mathrm{loop}}$. The proof follows the same line as in the case $e_0 \notin E_{\mathrm{loop}}$. Let $e_0 = \{v\}$ be incident to $v$. We prove only the case $v \neq v_0$. Defining $b$ as above, (48) is valid with

$$\begin{aligned} \frac{Z_{v_0,b}}{Z_{v_0,a}} &= \frac{\Gamma((a_v+1)/2)\Gamma((a_{e_0}+1)/2)\Gamma(a_{e_0}+2)}{4\Gamma((a_v+3)/2)\Gamma((a_{e_0}+3)/2)\Gamma(a_{e_0})} \\ &= \frac{a_{e_0}(a_{e_0}+1)}{(a_v+1)(a_{e_0}+1)} = \frac{a_{e_0}}{a_v+1}; \end{aligned} \tag{50}$$

here we used again the identity $\Gamma(z+1) = z\Gamma(z)$. The claim follows. □

Recall the definitions (5) and (6) of $k_v(\pi)$ and $k_e(\pi)$ for a finite admissible path $\pi$ in $G$. Abbreviate $k_v := k_v(\pi)$, $k_e := k_e(\pi)$. For $x = (x_e)_{e\in E} \in \Delta$, denote by $Q_x(\pi)$ the probability that the reversible Markov chain with transition probabilities induced by the weights $(x_e)_{e\in E}$ on the edges traverses the path $\pi$. Note that if $\pi$ is a *closed* path, that is, if the starting point and endpoint of $\pi$ agree, then $Q_x(\pi)$ is independent of the starting point of $\pi$. An argument as in the proof of Proposition 4.6 yields the following:

PROPOSITION 4.7. *For any finite admissible path $\pi$ starting at $v_0$, we have*

$$\begin{aligned}\int_\Delta Q_x(\pi)\phi_{v_0,a}(x)\,d\sigma(x)\\ = \frac{[\prod_{e\in E\setminus E_{\text{loop}}}\prod_{i=0}^{k_e-1}(a_e+i)][\prod_{e\in E_{\text{loop}}}\prod_{i=0}^{k_e/2-1}(a_e+2i)]}{\prod_{i=0}^{k_{v_0}-1}(a_{v_0}+2i)\prod_{v\in V\setminus\{v_0\}}\prod_{i=0}^{k_v-1}(a_v+1+2i)}.\end{aligned} \tag{51}$$

*For any finite admissible path $\pi$ with the same starting point and endpoint which avoids $v_0$, we have*

$$\begin{aligned}\int_\Delta Q_x(\pi)\phi_{v_0,a}(x)\,d\sigma(x)\\ = \frac{[\prod_{e\in E\setminus E_{\text{loop}}}\prod_{i=0}^{k_e-1}(a_e+i)][\prod_{e\in E_{\text{loop}}}\prod_{i=0}^{k_e/2-1}(a_e+2i)]}{\prod_{v\in V}\prod_{i=0}^{k_v-1}(a_v+1+2i)}.\end{aligned} \tag{52}$$

*Here the empty product is defined to be* 1.

If $\pi$ is a closed path, we call $Q_x(\pi)$ a *cycle probability*. The transition probabilities of a Markov chain with finite state space $V$ that visits every state with probability 1 are completely determined by all its cycle probabilities (see, e.g., [7], Corollary on page 116).

4.5. *Simulating from the posterior.* In this subsection we show how the posterior distribution of the unknown stationary distribution for the underlying Markov chain can be simulated using reinforced random walks.

Suppose our posterior distribution is $\mathbb{P}_{v_0,a} = \phi_{v_0,a}\,d\sigma$. Let $X^{(i)} := (X_n^{(i)})_{n\geq 0}$, $i \geq 1$, be independent reinforced random walks with the same initial edge weights $a = (a_e)_{e\in E}$. Let $Z_n^{(i)} := (X_0^{(i)}, X_1^{(i)}, \dots, X_n^{(i)})$ and recall that $k_e(Z_n^{(i)})$ equals the number of traversals of edge $e$ by the process $X^{(i)}$ up to time $n$.

PROPOSITION 4.8. *For any interval $I \subseteq \mathbb{R}$ and all $e \in E$, we have*

$$\lim_{n\to\infty}\lim_{m\to\infty}\frac{1}{m}\left|\left\{i\leq m : \frac{k_e(Z_n^{(i)})}{n}\in I\right\}\right| = \mathbb{P}_{v_0,a}(x_e\in I) \qquad a.s. \tag{53}$$

TABLE 1
*Degrees of freedom for independent, reversible and full Markov specification*

| $\|V\|$ | 3 | 4 | 5 | 10 | 20 | 50 | 100 | 1000 |
|---|---|---|---|---|---|---|---|---|
| Independent $\|V\|-1$ | 2 | 3 | 4 | 9 | 19 | 49 | 99 | 999 |
| Reversible $\|V\|(\|V\|-1)/2-1$ | 2 | 5 | 9 | 44 | 189 | 1224 | 4949 | 499499 |
| Full Markov $\|V\|(\|V\|-1)$ | 6 | 12 | 20 | 90 | 380 | 2450 | 9900 | 999000 |

PROOF. For every $n$, the random variables $k_e(Z_n^{(i)})/n$, $i \geq 1$, are i.i.d. Hence, by the Glivenko–Cantelli theorem, a.s. for all $x \in \mathbb{R}$,

$$\lim_{m\to\infty} \frac{1}{m}\left|\left\{i \leq m : \frac{k_e(Z_n^{(i)})}{n} \leq x\right\}\right| = P_{v_0,a}\left(\frac{k_e(Z_n)}{n} \leq x\right) = \mathbb{Q}_{v_0,a}\left(\frac{k_e(Z_n)}{n} \leq x\right). \tag{54}$$

For the last equality we used (19). Since $\mathbb{Q}_{v_0,a}$ is a mixture of Markov chains, $k_e(Z_n)/n$ converges to the normalized weight of the edge $e$ $\mathbb{Q}_{v_0,a}$-a.s. and, hence, weakly. Since the limiting distribution is continuous,

$$\lim_{n\to\infty} \mathbb{Q}_{v_0,a}\left(\frac{k_e(Z_n)}{n} \leq x\right) = \mathbb{P}_{v_0,a}(x_e \leq x), \tag{55}$$

and the claim follows. □

PROPOSITION 4.9. *For all $e \in E$,*

$$\lim_{n\to\infty} \int \frac{k_e(Z_n)}{n}\, dP_{v_0,a} = \int x_e \, d\mathbb{P}_{v_0,a}. \tag{56}$$

PROOF. By (19), $P_{v_0,a} = \mathbb{Q}_{v_0,a}$. Since the proportion $k_n(e)/n$ converges $\mathbb{Q}_{v_0,a}$-a.s. to the normalized weight of the edge $e$, the claim follows from the dominated convergence theorem. □

REMARK 4.10. The Markov chain with distribution induced by the edge weights $(x_e)_{e\in E} \in \Delta$ has the stationary distribution $\nu(v) = \frac{x_v}{2} = \frac{1}{2}\sum_{e\in E_v} x_e$. Thus, Propositions 4.8 and 4.9 allow simulation of the $\mathbb{P}_{v_0,a}$-distribution and the mean of $\nu(v)$.

**5. Applications.** Reversibility can serve as a natural intermediate between independence and fully nonparametric Markovian dependence. On

TABLE 2
*The humane HLA-B gene. Part of the DNA sequence of length* 3370

| | | | | | | |
|---|---|---|---|---|---|---|
| 1 | tggtgtagga | gaagagggat | caggacgaag | tcccaggtcc | cggacggggc | tctcagggtc |
| 61 | tcaggctccg | agggccgcgt | ctgcaatggg | gaggcgcagc | gttggggatt | ccccactccc |
| 121 | ctgagtttca | cttcttctcc | caacttgtgt | cgggtccttc | ttccaggata | ctcgtgacgc |
| 181 | gtccccactt | cccactccca | ttgggtattg | gatatctaga | gaagccaatc | agcgtcgccg |
| 241 | cggtcccagt | tctaaagtcc | ccacgcaccc | acccggactc | agagtctcct | cagacgccga |
| 301 | gatgctggtc | atggcgcccc | gaaccgtcct | cctgctgctc | tcggcggccc | tggccctgac |
| 361 | cgagacctgg | gccggtgagt | gcgggtcggg | agggaaatgg | cctctgccgg | gaggagcgag |
| 421 | gggaccgcag | gcgggggcgc | aggacctgag | gagccgcgcc | gggaggaggg | tcgggcgggt |
| 481 | ctcagcccct | cctcaccccc | aggctcccac | tccatgaggt | atttctacac | ctccgtgtcc |
| 541 | cggcccggcc | gcggggagcc | ccgcttcatc | tcagtgggct | acgtggacga | cacccagttc |
| 601 | gtgaggttcg | acagcgacgc | cgcgagtccg | agagaggagc | cgcgggcgcc | gtggatagag |
| 661 | caggaggggc | cggagtattg | ggaccggaac | acacagatct | acaaggccca | ggcacagact |
| 721 | gaccgagaga | gcctgcggaa | cctgcgcggc | tactacaacc | agagcgaggc | cggtgagtga |
| 781 | ccccggcccg | gggcgcaggt | cacgactccc | catcccccac | gtacggcccg | ggtcgccccg |
| 841 | agtctccggg | tccgagatcc | gcctccctga | ggccgcggga | cccgcccaga | ccctcgaccg |
| 901 | gcgagagccc | caggcgcgtt | tacccggttt | cattttcagt | tgaggccaaa | atccccgcgg |
| 961 | gttggtcggg | gcggggcggg | gctcggggga | ctgggctgac | cgcggggccg | gggccagggt |
| 1021 | ctcacaccct | ccagagcatg | tacggctgcg | acgtggggcc | ggacgggcgc | ctcctccgcg |
| 1081 | ggcatgacca | gtacgcctac | gacggcaagg | attacatcgc | cctgaacgag | gacctgcgct |
| 1141 | cctggaccgc | cgcggacacg | gcggctcaga | tcacccagcg | caagtgggag | gcggcccgtg |
| 1201 | aggcggagca | gcggagagcc | tacctggagg | gcgagtgcgt | ggagtggctc | cgcagatacc |
| 1261 | tggagaacgg | gaaggacaag | ctggagcgcg | ctggtaccag | gggcagtggg | gagccttccc |
| 1321 | catctcctat | aggtcgccgg | ggatggcctc | ccacgagaag | aggaggaaaa | tgggatcagc |
| 1381 | gctagaatgt | cgccctccgt | tgaatggaga | atggcatgag | ttttcctgag | tttcctctga |
| 1441 | gggccccctc | ttctctctag | acaattaagg | aatgacgtct | ctgaggaaat | ggaggggaag |
| 1501 | acagtcccta | gaatactgat | caggggtccc | ctttgacccc | tgcagcagcc | ttgggaaccg |
| 1561 | tgacttttcc | tctcaggcct | tgttctctgc | ctcacactca | gtgtgtttgg | ggctctgatt |
| 1621 | ccagcacttc | tgagtcactt | tacctccact | cagatcagga | gcagaagtcc | ctgttccccg |
| 1681 | ctcagagact | cgaactttcc | aatgaatagg | agattatccc | aggtgcctgc | gtccaggctg |
| 1741 | gtgtctgggt | tctgtgcccc | ttccccaccc | caggtgtcct | gtccattctc | aggctggtca |
| 1801 | catgggtggt | cctagggtgt | cccatgaaag | atgcaaagcg | cctgaatttt | ctgactcttc |
| 1861 | ccatcagacc | ccccaaagac | acacgtgacc | caccacccca | tctctgacca | tgaggccacc |
| 1921 | ctgaggtgct | gggccctggg | tttctaccct | gcggagatca | cactgacctg | gcagcgggat |
| 1981 | ggcgaggacc | aaactcagga | cactgagctt | gtggagacca | gaccagcagg | agatagaacc |
| 2041 | ttccagaagt | gggcagctgt | ggtggtgcct | tctggagaag | agcagagata | cacatgccat |
| 2101 | gtacagcatg | aggggctgcc | gaagcccctc | accctgagat | ggggtaagga | gggggatgag |
| 2161 | gggtcatatc | tcttctcagg | gaaagcagga | gcccttcagc | agggtcaggg | cccctcatct |
| 2221 | tcccctcctt | tcccagagcc | gtcttcccag | tccaccgtcc | ccatcgtggg | cattgttgct |
| 2281 | ggcctggctg | tcctagcagt | tgtggtcatc | ggagctgtgg | tcgctgctgt | gatgtgtagg |
| 2341 | aggaagagtt | caggtaggga | aggggtgagg | ggtggggtct | gggttttctt | gtcccactgg |
| 2401 | gggtttcaag | ccccaggtag | aagtgttccc | tgcctcatta | ctgggaagca | gcatgcacac |
| 2461 | aggggctaac | gcagcctggg | accctgtgtg | ccagcactta | ctcttttgtg | cagcacatgt |
| 2521 | gacaatgaag | gatggatgta | tcaccttgat | ggttgtggtg | ttggggtcct | gattccagca |
| 2581 | ttcatgagtc | aggggaaggt | ccctgctaag | gacagacctt | aggagggcag | ttggtccagg |
| 2641 | acccacactt | gctttcctcg | tgtttcctga | tcctgccctg | ggtctgtagt | catacttctg |
| 2701 | gaaattcctt | ttgggtccaa | gactaggagg | ttcctctaag | atctcatggc | cctgcttcct |
| 2761 | cccagtgccc | tcacaggaca | ttttcttccc | acaggtggaa | aaggagggag | ctactctcag |
| 2821 | gctgcgtgta | agtggtgggg | gtgggagtgt | ggaggagctc | acccacccca | taattcctcc |
| 2881 | tgtcccacgt | ctcctgcggg | ctctgaccag | gtcctgtttt | tgttctactc | caggcagcga |
| 2941 | cagtgcccag | ggctctgatg | tgtctctcac | agcttgaaaa | ggtgagattc | ttggggtcta |
| 3001 | gagtgggtgg | ggtggcgggt | ctgggggtgg | gtggggcaga | ggggaaaggc | ctgggtaatg |
| 3061 | gggattcttt | gattgggatg | tttcgcgtgt | gtggtgggct | gtttagagtg | tcatcgctta |
| 3121 | ccatgactaa | ccagaatttg | ttcatgactg | ttgttttctg | tagcctgaga | cagctgtctt |
| 3181 | gtgagggact | gagatgcagg | atttcttcac | gcctcccctt | tgtgacttca | agagcctctg |
| 3241 | gcatctcttt | ctgcaaaggc | acctgaatgt | gtctgcgtcc | ctgttagcat | aatgtgagga |
| 3301 | ggtggagaga | cagcccaccc | ttgtgtccac | tgtgacccct | gttcgcatgc | tgacctgtgt |
| 3361 | ttcctcccca | | | | | |

$|V|$ states, with no restrictions the number of free parameters is $|V|-1$ with

TABLE 3
*Occurrences $N_{ij}$ of the string $ij$ for $i, j \in \{a, c, g, t\}$*

| | **a** | **c** | **g** | **t** |
|---|---|---|---|---|
| $a$ | 91 | 160 | 261 | 108 |
| $c$ | 213 | 351 | 161 | 249 |
| $g$ | 251 | 224 | 388 | 201 |
| $t$ | 66 | 239 | 254 | 152 |

independence, $|V|(|V|-1)$ for full Markov and $\frac{|V|(|V|-1)}{2}-1$ for reversibility. As Table 1 indicates, these numbers vary widely for $|V|$ large.

In this section we illustrate the use of our priors for testing a variety of simple hypotheses. Table 2 shows a genetic data set from the DNA sequence of the humane HLA-B gene. This gene plays a central role in the immune system. The data displayed in Table 2 is downloaded from the webpage of the National Center for Biotechnology Information (`www.ncbi.nlm.nih.gov/genome/guide/human/`).

In Example A, we test i.i.d. $\frac{1}{4}$ versus i.i.d. for the DNA-data. In Example B, we test i.i.d. versus reversible. In Example C, we test reversible versus full Markov. In Example D, we compare i.i.d. with full Markov.

Let $n_a$, $n_c$, $n_g$ and $n_t$ denote the number of occurrences of $a$, $c$, $g$ and $t$, respectively, in the data displayed in Table 2. Then

$$n_a = 621, \qquad n_c = 974, \qquad n_g = 1064, \qquad n_t = 711. \tag{57}$$

EXAMPLE A. A Bayes test of $H_0$: i.i.d.$(\frac{1}{4})$ versus $H_1$: i.i.d.(unknown). A "standard" test can be based on the Bayes factor

$$\frac{P(data|H_0)}{P(data|H_1)}.$$

See [9] for an extensive discussion. For $H_1$, we use a Dirichlet$(1,1,1,1)$ prior. This yields

$$P(data|H_0) = \left(\frac{1}{4}\right)^{3370} \approx 1.142429015368253 \cdot 10^{-2029},$$

$$\begin{aligned} P(data|H_1) &= \frac{\Gamma(4)\Gamma(n_a+1)\Gamma(n_c+1)\Gamma(n_g+1)\Gamma(n_t+1)}{\Gamma(n_a+n_c+n_g+n_t+4)} \\ &= \frac{\Gamma(4)\Gamma(622)\Gamma(975)\Gamma(1065)\Gamma(712)}{\Gamma(3374)} \\ &\approx 1.140417804695619 \cdot 10^{-1999}, \end{aligned}$$

$$\frac{P(data|H_0)}{P(data|H_1)} \approx 1.00176 \cdot 10^{-30}.$$

Thus, $H_0$ is strongly rejected. This is not surprising since the observed numbers of $a$, $c$, $g$, $t$ are $n_a = 621$, $n_c = 974$, $n_g = 1064$, $n_t = 711$, respectively.

EXAMPLE B. A Bayes test of $H_0$: i.i.d.(unknown) versus $H_1$: reversible. Here we use a Dirichlet$(1,1,1,1)$ prior for the null hypothesis and the prior based on the complete graph $K_4$ with loops (see Figure 5) and all edge weights equal to 1. The probability $P(data|H_0)$ is calculated in Example A. In order to calculate $P(data|H_1)$, we first determine the transition counts $k_e$ for our data (see Table 4) and also $k_v = n_v - \delta_a(v)$:

$$k_a = 620, \qquad k_c = 974, \qquad k_g = 1064, \qquad k_t = 711. \tag{58}$$

We abbreviate $E' = \{\{a,c\},\{a,g\},\{a,t\},\{c,g\},\{c,t\},\{g,t\}\}$. By the first part of Proposition 4.7,

$$\begin{aligned} P(data|H_1) &= \frac{\prod_{e\in E'}\prod_{i=0}^{k_e-1}(1+i)\prod_{j\in\{a,c,g,t\}}\prod_{i=0}^{k_{\{j\}}/2-1}(1+2i)}{\prod_{i=0}^{k_t-1}(4+2i)\prod_{j\in\{a,c,g\}}\prod_{i=0}^{k_j-1}(5+2i)} \\ &= (373)!(512)!(174)!(385)!(488)!(455)! \\ &\quad \times \prod_{i=0}^{90}(1+2i)\prod_{i=0}^{350}(1+2i)\prod_{i=0}^{387}(1+2i)\prod_{i=0}^{151}(1+2i) \\ &\quad \times \left(\prod_{i=0}^{710}(4+2i)\prod_{i=0}^{619}(5+2i)\prod_{i=0}^{973}(5+2i)\prod_{i=0}^{1063}(5+2i)\right)^{-1} \\ &\approx 2.166939224648291 \cdot 10^{-1961}. \end{aligned} \tag{59}$$

So the Bayes factor is

$$\frac{P(data|H_0)}{P(data|H_1)} \approx 5.2628 \cdot 10^{-39}$$

and the null hypothesis is strongly rejected.

TABLE 4
*The undirected transition counts $k_{\{i,j\}}$, $i,j \in \{a,c,g,t\}$*

| | **a** | **c** | **g** | **t** |
|---|---|---|---|---|
| $a$ | 182 | 373 | 512 | 174 |
| $c$ | 373 | 702 | 385 | 488 |
| $g$ | 512 | 385 | 776 | 455 |
| $t$ | 174 | 488 | 455 | 304 |

EXAMPLE C. A Bayes test of $H_0$: reversible versus $H_1$: full Markov.

Here we use our conjugate prior on reversible chains with all constants chosen as one. We use product Dirichlet measure for the rows in the full Markov case. This yields

$$\begin{aligned} P(data|H_1) &= \prod_{i\in\{a,c,g,t\}} \Gamma(4) \frac{\prod_{j\in\{a,c,g,t\}} \Gamma(N_{ij}+1)}{\Gamma(k_i+4)} \\ &= \Gamma(4)^4 \cdot \frac{\Gamma(92)\Gamma(161)\Gamma(262)\Gamma(109)}{\Gamma(624)} \\ &\quad \times \frac{\Gamma(214)\Gamma(352)\Gamma(162)\Gamma(250)}{\Gamma(978)} \\ &\quad \times \frac{\Gamma(252)\Gamma(225)\Gamma(389)\Gamma(202)}{\Gamma(1068)} \cdot \frac{\Gamma(67)\Gamma(240)\Gamma(255)\Gamma(153)}{\Gamma(715)} \\ &\approx 4.16382063735625 \cdot 10^{-1956}. \end{aligned}$$

The probability $P(data|H_0)$ was calculated in Example B. Hence,

$$\frac{P(data|H_0)}{P(data|H_1)} \approx 5.20421 \cdot 10^{-6}.$$

We see that a straightforward Bayes test rejects reversibility.

EXAMPLE D. A Bayes test of $H_0$: i.i.d.(unknown) versus $H_1$: full Markov.

Using the Bayes factors computed above, the null hypothesis is strongly rejected:

$$\frac{P(data|H_0)}{P(data|H_1)} \approx 2.73887 \cdot 10^{-44}.$$

Of course, an i.i.d. process is a reversible Markov chain.

In using the Dirichlet prior for testing uniformity with multinomial data and for testing independence in contingency tables, I. J. Good found the symmetric Dirichlet prior with density proportional to $\prod_{i=1}^{d} x_i^{c-1}$ an important tool. Good's many insights into these testing problems may be accessed through his book [9] and the survey article [10].

We have used the analog of the symmetric Dirichlet for the reversible Markov chain context with all edge weights $a_e$ equal to a constant $c$ say. As $c$ tends to infinity, this prior tends to a point mass supported on the simple random walk on the graph. As $c$ tends to zero, this prior tends to an improper prior which gives the maximum likelihood as its posterior.

Good also worked with $c$-mixtures of symmetric Dirichlet priors. We suspect that parallel, useful things can be done in our case as well.

We have not found *any* literature about statistical analysis of reversible Markov chains with unknown transitions and append two data analytic remarks here. First, under reversibility, the count $N_{vv'}$ of $v$ to $v'$ transitions has the same expectation as the count $N_{v'v}$ of $v'$ to $v$ transitions, namely, $\nu(v)k(v,v')$. This suggests looking at ratios $N_{vv'}/N_{v'v}$ or differences $N_{vv'} - N_{v'v}$. For example, from Table 3, $N_{ac}/N_{ca} = 160/213$, $N_{ag}/N_{ga} = 261/251$, $N_{at}/N_{ta} = 108/66$, $N_{cg}/N_{gc} = 161/224$, $N_{ct}/N_{tc} = 249/239$, $N_{gt}/N_{tg} = 201/254$; most of these are way off.

In large samples, these counts have limiting normal distributions by results of Höglund [11]. A second data analytic tool would be to estimate the stationary distribution [perhaps by the method of moments estimator $\hat{\nu}(v) = \frac{1}{n}|\{i \le n : X_i = v\}|$] and also estimate the transition matrix, and then compare $\hat{\nu}(v)\hat{k}(v,v')$ with $\hat{\nu}(v')\hat{k}(v',v)$.

An interesting problem not tackled here is finding natural priors on the set of reversible Markov chains *with a fixed stationary distribution*. For definiteness, consider the uniform stationary distribution. Then the problem is to put a prior on $\mathcal{S}(n)$, the symmetric doubly stochastic $n \times n$ matrices. We make two remarks. First, determining the Euclidean volume of $\mathcal{S}(n)$ is a long-standing open problem; see [1] for recent results. Second, $\mathcal{S}(n)$ is a compact, convex subset of $\mathbb{R}^{n^2}$. Its extreme points are well known to be the *symmetrized* permutation matrices (see [17]). Thus, if $\pi$ is a permutation matrix on $n$ letters with $e(\pi)$ the usual $n \times n$ permutation matrix, let $\tilde{e}(\pi) = \frac{1}{2}[e(\pi) + e(\pi^{-1})]$. The extreme points of $\mathcal{S}(n)$ are $(\tilde{e}(\pi))$ as $\pi$ ranges over permutations in $S_n$. We may put a prior on $\mathcal{S}(n)$ by taking a random convex combination of the $\tilde{e}(\pi)$. Alas, $\mathcal{S}(n)$ is *not* a simplex, so symmetric weights on the extreme points may not lead to symmetric measures on $\mathcal{S}(n)$.

**Acknowledgment.** We would like to thank Franz Merkl for some interesting discussions.

## REFERENCES

[1] Chan, C. S., Robbins, D. P. and Yuen, D. S. (2000). On the volume of a certain polytope. *Experiment. Math.* **9** 91–99. MR1758803

[2] Dalal, S. R. and Hall, W. J. (1983). Approximating priors by mixtures of natural conjugate priors. *J. Roy. Statist. Soc. Ser. B* **45** 278–286. MR0721753

[3] Diaconis, P. (1988). Recent progress on de Finetti's notions of exchangeability. In *Bayesian Statistics 3* (J.-M. Bernardo, M. H. DeGroot, D. V. Lindley and A. F. M. Smith, eds.) 111–125. Oxford Univ. Press, New York. MR1008047

[4] Diaconis, P. and Freedman, D. (1980). de Finetti's theorem for Markov chains. *Ann. Probab.* **8** 115–130. MR0556418

[5] Diaconis, P. and Ylvisaker, D. (1985). Quantifying prior opinion (with discussion). In *Bayesian Statistics 2* (J.-M. Bernardo, M. H. DeGroot, D. V. Lindley and A. F. M. Smith, eds.) 133–156. North-Holland, Amsterdam. MR0862488

[6] Fortini, S., Ladelli, L., Petris, G. and Regazzini, E. (2002). On mixtures of distributions of Markov chains. *Stochastic Process. Appl.* **100** 147–165. MR1919611

[7] Freedman, D. A. (1962). Mixtures of Markov processes. *Ann. Math. Statist.* **33** 114–118. MR0137156

[8] Giblin, P. J. (1981). *Graphs, Surfaces and Homology. An Introduction to Algebraic Topology*, 2nd ed. Chapman and Hall, London. MR0643363

[9] Good, I. J. (1968). *The Estimation of Probabilities. An Essay on Modern Bayesian Methods.* MIT Press, Cambridge, MA. MR0245135

[10] Good, I. J. and Crook, J. F. (1987). The robustness and sensitivity of the mixed-Dirichlet Bayesian test for "independence" in contingency tables. *Ann. Statist.* **15** 670–693. MR0888433

[11] Höglund, T. (1974). Central limit theorems and statistical inference for finite Markov chains. *Z. Wahrsch. Verw. Gebiete* **29** 123–151. MR0373201

[12] Keane, M. S. (1990). Solution to problem 288. *Statist. Neerlandica* **44** 95–100.

[13] Keane, M. S. and Rolles, S. W. W. (2000). Edge-reinforced random walk on finite graphs. In *Infinite Dimensional Stochastic Analysis* (Ph. Clément, F. den Hollander, J. van Neerven and B. de Pagter, eds.) 217–234. R. Neth. Acad. Arts Sci., Amsterdam. MR1832379

[14] Maurer, S. B. (1976). Matrix generalizations of some theorems on trees, cycles and cocycles in graphs. *SIAM J. Appl. Math.* **30** 143–148. MR0392635

[15] Pemantle, R. (1988). Phase transition in reinforced random walk and RWRE on trees. *Ann. Probab.* **16** 1229–1241. MR0942765

[16] Rolles, S. W. W. (2003). How edge-reinforced random walk arises naturally. *Probab. Theory Related Fields* **126** 243–260. MR1990056

[17] Stanley, R. P. (1978). Generating functions. In *Studies in Combinatorics* (G.-C. Rota, ed.) 100–141. Math. Assoc. Amer., Washington. MR0513004

[18] Zabell, S. L. (1982). W. E. Johnson's "sufficientness" postulate. *Ann. Statist.* **10** 1090–1099. MR0673645

[19] Zabell, S. L. (1995). Characterizing Markov exchangeable sequences. *J. Theoret. Probab.* **8** 175–178. MR1308676

Department of Statistics
Stanford University
Stanford, California 94305-4065
USA
E-mail: diaconis@math.stanford.edu

Eindhoven University of Technology
Department of Mathematics
and Computer Science
P.O. Box 513
5600 MB Eindhoven
The Netherlands
E-mail: srolles@win.tue.nl